\documentclass[11pt]{article}

\usepackage[]{graphicx,float,latexsym,times}
\usepackage{amsfonts,amstext,amsmath,amssymb,amsthm,bm}
\usepackage{caption2}

\long\def\symbolfootnote[#1]#2{\begingroup%
\def\thefootnote{\fnsymbol{footnote}}\footnote[#1]{#2}\endgroup}

\usepackage{titlesec} 

\titleformat{\section}{\large\bfseries}{\thesection.}{.5em}{}
\titlespacing*{\section}{0pt}{*3}{*2}
\titleformat{\subsection}{\normalfont\bfseries}{\thesubsection.}{.5em}{}
\titlespacing*{\subsection} {0pt}{*3}{*2}
\titleformat{\subsubsection}{\normalfont\bfseries}{\thesubsubsection.}{.5em}{}
\titlespacing*{\subsubsection} {0pt}{*3}{*2}

\theoremstyle{plain} 
\newtheorem{theorem}{Theorem}[section]

\newtheorem{corollary}{Corollary}[section]

\theoremstyle{definition} 

\newcommand{\To}{\rightarrow}

\newcommand{\qmq}[1]{\quad\mbox{#1}\quad}

\numberwithin{equation}{section} 

\usepackage[authoryear]{natbib}

\begin{document}

\title{\textbf{\Large A Proof of the Bomber Problem's Spend-It-All Conjecture}}

\date{}

\maketitle


\author{
\begin{center}
\vskip -1cm

\textbf{\large Jay Bartroff}

Department of Mathematics, University of Southern California, Los Angeles, California, USA
\end{center}
}

\symbolfootnote[0]{\normalsize Address correspondence to Jay Bartroff,
Department of Mathematics, University of Southern California, KAP 108, Los Angeles,
CA 90089, USA; E-mail: bartroff@usc.edu}

{\small \noindent\textbf{Abstract:} The Bomber Problem concerns optimal sequential allocation of partially effective ammunition~$x$ while under attack from enemies arriving according to a Poisson process over a time interval of length~$t$. In the doubly-continuous setting, in certain regions of $(x,t)$-space we are able to solve the integral equation defining the optimal survival probability and find  the optimal allocation function~$K(x,t)$ exactly in these regions. As a consequence, we complete the proof of the ``spend-it-all'' conjecture of \citet{Bartroff10} which gives the boundary of the region where $K(x,t)=x$.}

\section{INTRODUCTION}
The Bomber Problem, first posed by \citet{Klinger68}, concerns an aircraft equipped with $x$ units
of ammunition that is $t$ time units away from its final destination. It is confronted by enemy airplanes whose appearance
is driven by a time-homogenous Poisson process with known intensity, taken to be 1. In the encounters, which are assumed to be instantaneous, if the Bomber fires $y$ of its available~$x$ units of ammunition, then the probability that the Bomber survives the encounter is given by
\begin{equation}\label{a-function}
a(y)=1-(1-u)e^{-y}
\end{equation}
for some fixed $0\le u< 1$. This can be interpreted as the Bomber's
$y$ units of ammunition destroying the enemy with
probability~$1-e^{-y}$, while otherwise allowing the enemy to launch a counterattack
which succeeds with probability~$1-u$. The optimal amount of ammunition that should be spent when the Bomber is confronted by an enemy while in
``state''~$(x,t)$ in order to maximize the probability of reaching its destination is denoted $K(x,t)$, and the optimal probability is denoted $P(x,t)$. Of central interest in the Bomber Problem is a set of conjectures concerning monotonicity of $K(x,t)$; see \citet{Bartroff10c} and references therein for a description of the conjectures and their statuses.

In this paper we adopt the doubly-continuous setting where $x$ and $t$ are both assumed to be continuous variables, and show in Theorem~\ref{thm:KP} that $K(x,t)$ and $P(x,t)$ can be solved for exactly in certain regions of $(x,t)$-space. This is done by solving the integral equation
\begin{equation}\label{int_eq}
P(x,t)=e^{-t}\left(1+\int_0^t\max_{0\le y\le x}a(y)P(x-y,s)e^sds\right),
\end{equation}
of which \citet[][Corollary~2.2]{Bartroff10c} showed that the the optimal survival probability~$P(x,t)$ is the unique solution, and showing that the maximum in the integrand of (\ref{int_eq}) is uniquely achieved, giving $K(x,t)$. This allows us in Corollary~\ref{cor:SIA} to complete the proof of a conjecture by \citet{Bartroff10} that the ``spend-it-all'' region, i.e., the set of all $(x,t)$ such that $K(x,t)=x$, is given by
\begin{equation}\label{R1}
R_1=\{(x,t)\in\mathbb{R}^+\times \mathbb{R}^+: x\le f_u(t)\},
\end{equation} where
\begin{equation}\label{f}
f_u(t) =\begin{cases}
\log(1+u/(e^{tu}-1)),&0<u\le 1\\
\log(t^{-1}+1),&u=0.
\end{cases}
\end{equation}
Regarding the Bomber Problem's monotonicity conjectures we only remark that the form of $K(x,t)$ found in (\ref{K}) is increasing in $x$, and hence does not violate the conjecture that $K(x,t)$ obeys this everywhere, although this remains unproved in general at this time; see \citet{Bartroff10c}.

\section{Results}

\begin{theorem}\label{thm:KP}
Let $R_1$ be as in (\ref{R1}) and
$$R_2=\{(x,t)\in\mathbb{R}^+\times \mathbb{R}^+: f_u(t)<x\le 2f_u(t)\}.$$
For all $0\le u<1$,
\begin{equation}\label{K}
K(x,t)=\begin{cases}
x,&(x,t)\in R_1\\
(x+f_u(t))/2, &(x,t)\in R_2.
\end{cases} 
\end{equation} For all $0<u< 1$,
\begin{equation}\label{P}
P(x,t)=e^{-t}\times\begin{cases}
1+a(x)(e^{tu}-1)/u,&(x,t)\in R_1\\
1+\frac{a(x)}{e^x-1}+\int_{f_u(x/u)/u}^t u^{-1}\left(\sqrt{e^{s u}-1+u}-(1-u)e^{-x/2}\sqrt{e^{s u}-1}\right)^2 ds , &(x,t)\in R_2,
\end{cases} 
\end{equation} and (\ref{P}) holds for $u=0$ by taking the limit of the right-hand-side as $u\To 0$.
\end{theorem}

\noindent\textit{Proof.}  We show that the right-hand-side of (\ref{P}) satisfies the integral equation~(\ref{int_eq}), and then show that the maximum in (\ref{int_eq}) is uniquely achieved at (\ref{K}). Assume that $0<u< 1$. The proof for $u=0$ is exactly the same after replacing all quantities involving $u$ by their limit as $u\To 0$. For example, the integrand in the second case of (\ref{P}) becomes
$$\left(\sqrt{s+1}-e^{-x/2}\sqrt{s}\right)^2.$$
Letting $\overline{P}(x,t)=e^t P(x,t)$, (\ref{int_eq}) becomes
\begin{equation}\label{int_eq2}
\overline{P}(x,t)=1+\int_0^t \max_{0\le y\le x}a(y)\overline{P}(x-y,s) ds.
\end{equation} A simple but useful fact is that, for fixed $x,t,B>0$, the function 
\begin{equation}\label{uni}
y\mapsto a(y)(1+Ba(x-y))\qmq{is unimodal about}y^*=\frac{x+\log(1+1/B)}{2},\end{equation} which can be verified by basic calculus. 

Fix $(x,t)\in R_1$ and let $0\le s \le t$ and $$G_1(y,s)=a(y)\left[1+a(x-y)(e^{s u}-1)/u\right].$$  In order to compute $\max_{0\le y\le x}G_1(y,s)$ we apply (\ref{uni}) with $B=(e^{s u}-1)/u$ and, using that $f_u$ is decreasing, we have $$\left.\frac{x+\log(1+1/B)}{2}\right|_{B=(e^{s u}-1)/u}=\frac{x+f_u(s )}{2}\ge \frac{x+f_u(t)}{2}\ge \frac{x+x}{2}=x.$$ Thus 
\begin{equation}\label{maxR1}
\max_{0\le y\le x}G_1(y,s)=G_1(x,s)=a(x)e^{s u},\end{equation} hence
$$1+\int_0^t \max_{0\le y\le x}G_1(y,s)ds =1+ \int_0^t a(x)e^{s u} ds =1+a(x)(e^{tu}-1)/u,$$ giving the first case of (\ref{P}), and also (\ref{K}) since the maximum is uniquely achieved at $y=x$.

Now fix $(x,t)\in R_2$. Let $v=1-u$ and 
\begin{align}
q(y,s)&=\left(\sqrt{e^{s u}-v}-ve^{-y/2}\sqrt{e^{s u}-1}\right)^2/u,\\
Q_2(y,s)&=1+\frac{a(y)}{e^y-1}+\int_{f_u(y/u)/u}^s q(y,r) dr ,\label{Q1}\\
Q(y,s)&=\begin{cases}
1+a(y)(e^{su}-1)/u,& (y,s)\in R_1\\
Q_2(y,s),& (y,s)\in R_2,
\end{cases}\\
G(s)&=\max_{0\le y\le x}a(y)Q(x-y,s).
\end{align} Let $0\le s\le t$ and note that $f_u^{-1}(y)=f_u(y/u)/u$. If $s\le f_u(x/u)/u$, then $(x-y,s)\in R_1$ for all $0\le y\le x$ since $x-y\le x\le f_u(s)$, hence
\begin{equation}
G(s)=a(x)e^{s u}\qmq{for} 0\le s\le f_u(x/u)/u\label{Gsmalls}
\end{equation} by (\ref{maxR1}). Now assume that $f_u(x/u)/u\le s\le t$. Then $$(x-y,s)\in \begin{cases}
R_1,& x-f_u(s)\le y\le x\\
R_2,&0\le y\le x-f_u(s),
\end{cases}$$ hence
\begin{equation}\label{Gtau2}
G(s)=\max\left\{\max_{0\le y\le x-f_u(s)} a(y)Q_2(x-y,s), \max_{x-f_u(s)\le y\le x} G_1(y,s)\right\}.
\end{equation} To compute the second term in (\ref{Gtau2}),  we again apply (\ref{uni}) with $B=(e^{s u}-1)/u$ but this time note that
\begin{equation}
\left.\frac{x+\log(1+1/B)}{2}\right|_{B=(e^{s u}-1)/u}=\frac{x+f_u(s)}{2}\in [x-f_u(s),x],\label{y*}\end{equation} hence
\begin{equation}
\max_{x-f_u(s)\le y\le x} G_1(y,s) =G_1((x+f_u(s ))/2,s) =q(x,s),\label{maxbigy}
\end{equation} after some simplification. Letting $G_2(y,s)=a(y)Q_2(x-y,s)$ and $G_2'(y,s)=(\partial/\partial y)G_2(y,s)$, we show that $G(s)$ is in fact equal to (\ref{maxbigy}) by showing that $G_2'(y,s)>0$ for all $0\le y\le x-f_u(s)$. 
Letting $Q_2'(y,s)=(\partial/\partial y)Q_2(y,s)$, using the fundamental theorem of calculus we have
\begin{align}
Q_2'(y,s)&=\frac{2v-ve^{-y}-e^y}{(e^y-1)^2}-\left(\frac{\partial}{\partial y} f_u(y/u)/u\right)q(y,f_u(y/u)/u)+\int_{f_u(y/u)/u}^s \frac{\partial}{\partial y}q(y,r)dr\nonumber\\
&=\frac{ve^{-y}}{e^y-1} +\frac{ve^{-y/2}}{u}\int_{f_u(y/u)/u}^s \sqrt{e^{r u}-1} \left(\sqrt{e^{r u}-v}-ve^{-y/2}\sqrt{e^{r u}-1}\right) dr.\label{Q1'}
\end{align} Letting $I_1(y,s)$ be the second term in (\ref{Q1'}) and $I_2(y,s)$ the integral in (\ref{Q1}), we have
\begin{align}\label{G'}
G_2'(y,s)&=a'(y)Q_2(x-y,s)-a(y)Q_2'(x-y,s)\nonumber\\
&=ve^{-y}\left(1+\frac{a(x-y)}{e^{x-y}-1}+I_2(x-y,s)\right)-a(y)\left(\frac{ve^{-x+y}}{e^{x-y}-1}+I_1(x-y,s)\right).
\end{align} 
The terms in (\ref{G'}) not involving integrals are
\begin{equation}\label{G1}
ve^{-y}\left(1+\frac{a(x-y)}{e^{x-y}-1}\right)-a(y)\left(\frac{ve^{-x+y}}{e^{x-y}-1}\right)=\frac{v(e^{x-2y}-e^{-x+y})}{e^{x-y}-1}>0
\end{equation} since, comparing exponents in the numerator and using that $y\le x-f_u(s)$, 
\begin{align*}
x-2y&=2x-3y-x+y \\
&\ge2x-3(x-f_u(t))-x+y\\
&=(3f_u(t)-x)-x+y\\
&> -x+y.
\end{align*} The terms in (\ref{G'}) involving integrals are
\begin{multline}
ve^{-y}I_2(x-y,s)-a(y)I_1(x-y,s)=\\
\frac{ve^{-(x-y)/2}}{u}\int_{f_u((x-y)/u)/u}^s \left(\sqrt{e^{r u}-v}-ve^{-(x-y)/2}\sqrt{e^{r u}-1}\right)\sqrt{e^{r u}-1}\left(e^{x/2-3y/2+f_u(r)/2}-1\right)dr\ge 0,\label{G2}
\end{multline} since the integrand is nonnegative: The first two factors are clearly so and the third factor is as well since, examining the exponent and using that $y\le x-f_u(s)$ and $f_u$ is decreasing,
\begin{align*}
x/2-3y/2+f_u(r)/2&\ge x/2-3(x-f_u(s))/2+f_u(s)/2\\
&=-x+2f_u(s)\\
&\ge -x+2f_u(t)\\
&\ge 0
\end{align*} since $(x,t)\in R_2$. Combining (\ref{G1}) and (\ref{G2}) shows that $G_2'(y,s)>0$ for all $0\le y\le x-f_u(s)$, and hence that $G(s)$ is equal to (\ref{maxbigy}) for $f_u(x/u)/u\le s\le t$. Note that in this case the maximum in $G(s)$ is uniquely attained at $y$ equal to (\ref{y*}). Then, using (\ref{Gsmalls}), 
\begin{align*}
1+\int_0^t \max_{0\le y\le x}a(y)Q(x-y,s)ds&=1+\int_0^t G(s)ds\\
&=1+\int_0^{f_u(x/u)/u}a(x)e^{su}ds+ \int_{f_u(x/u)/u}^t q(x,s)ds\\
&=1+\frac{a(x)}{e^x-1}+ \int_{f_u(x/u)/u}^t q(x,s)ds\\
&=Q_2(x,t)\\
&=Q(x,t),
\end{align*} this last since $(x,t)\in R_2$, showing that $Q$ satisfies the integral equation in $R_2$ and hence showing that the second case of (\ref{K}) and (\ref{P}) hold.\qed

\bigskip

As a consequence of Theorem~\ref{thm:KP}, we complete the proof of the spend-it-all conjecture:
\begin{corollary} \label{cor:SIA}$K(x,t)=x$ if and only if $x\le f_u(t)$.
\end{corollary}

\noindent\textit{Proof.}  \citet[][Theorem~2.1]{Bartroff10} showed that $K(x,t)<x$ if $x>f_u(t)$, and Theorem~\ref{thm:KP} provides the converse.\qed

\section*{ACKNOWLEDGEMENTS}

This work was supported in part by grant DMS-0907241 from the National Science Foundation and a Faculty-in-Residence grant from the Albert and Elaine Borchard Foundation.



\def\cprime{$'$}

\end{document}